\documentclass{article}
\usepackage{amsmath,amsthm}
\usepackage{amssymb}
\usepackage[all]{xy}

\newcommand{\dpullback}[1][d]{\save*!/#1-1.8pc/#1:(-1,1)@^{|-}\restore}

\usepackage{url}
\usepackage{texdraw}
\input{txdtools}

\def\fatness{1.2}

\everytexdraw{%
	\drawdim pt \linewd 0.5 \textref h:C v:C
  \arrowheadsize l:4 w:3  
  \arrowheadtype t:V
}

\newcommand{\multFA}{%
	\bsegment
		\realmult {\fatness} {2.5} \xrad
		\realmult {\fatness} {6.25} \yrad
		\move (0 0) \lellip rx:{\xrad} ry:{\yrad}
		\move (0 30) \lellip rx:{\xrad} ry:{\yrad}
		\move (40 15) \lellip rx:{\xrad} ry:{\yrad}
		\realadd {15} { \yrad} \outB
		\realadd {15} {-\yrad} \outA
		\realadd {30} { \yrad} \inD
		\realadd {30} {-\yrad} \inC
		\realadd { 0} { \yrad} \inB
		\realadd { 0} {-\yrad} \inA
		\move (0 {\inD}) \clvec (25 {\inD})(20 {\outB})(40 {\outB}) 
		\move (0 {\inB}) \clvec (18 {\inB})(18 {\inC})(0 {\inC}) 
		\move (0 {\inA}) \clvec (25 {\inA})(20 {\outA})(40 {\outA}) 
		\savepos (40 15)(*ex *ey)
	\esegment
	\move (*ex *ey)
}
\newcommand{\unitFA}{%
	\bsegment
		\realmult {\fatness} {2.5} \xrad
		\realmult {\fatness} {6.25} \yrad
		\move (0 0)
		\move (-20 0)
		\move (0 0) \lellip rx:{\xrad} ry:{\yrad}
		\realadd { 0} { \yrad} \inB
		\realadd { 0} {-\yrad} \inA
		\move (0 {\inB}) \rlvec (-8 0) 
		\move (0 {\inA}) \rlvec (-8 0) 
	  \move (-8 {\inA}) \clvec (-20 {\inA})(-20 {\inB})(-8 {\inB}) 

		\savepos (0 0)(*ex *ey)
	\esegment
	\move (*ex *ey)
}
\newcommand{\comultFA}{%
	\bsegment
		\realmult {\fatness} {2.5} \xrad
		\realmult {\fatness} {6.25} \yrad
		\move (0 0) \lellip rx:{\xrad} ry:{\yrad}
		\move (40 -15) \lellip rx:{\xrad} ry:{\yrad}
		\move (40 15) \lellip rx:{\xrad} ry:{\yrad}
		\realadd { 0} { \yrad} \outB
		\realadd { 0} {-\yrad} \outA
		\realadd {15} { \yrad} \inD
		\realadd {15} {-\yrad} \inC
		\realadd {-15} { \yrad} \inB
		\realadd {-15} {-\yrad} \inA
		\move (40 {\inD}) \clvec (15 {\inD})(20 {\outB})(0 {\outB}) 
		\move (40 {\inB}) \clvec (22 {\inB})(22 {\inC})(40 {\inC}) 
		\move (40 {\inA}) \clvec (15 {\inA})(20 {\outA})(0 {\outA}) 
		\savepos (40 -15)(*ex *ey)
	\esegment
	\move (*ex *ey)
}
\newcommand{\counitFA}{%
	\bsegment
		\realmult {\fatness} {2.5} \xrad
		\realmult {\fatness} {6.25} \yrad
		\move (0 0)
		\move (20 0)
		\move (0 0) \lellip rx:{\xrad} ry:{\yrad}
		\realadd { 0} { \yrad} \inB
		\realadd { 0} {-\yrad} \inA
		\move (0 {\inB}) \rlvec (8 0) 
		\move (0 {\inA}) \rlvec (8 0) 
	  \move (8 {\inA}) \clvec (20 {\inA})(20 {\inB})(8 {\inB}) 

		\savepos (0 0)(*ex *ey)
	\esegment
	\move (*ex *ey)
}
\newcommand{\handleFA}{%
       \realmult {\fatness} {2.5} \xrad
        \realmult {\fatness} {6.25} \yrad
	\bsegment
		\rmove (40 -15)
    \bsegment
         \move (-40 15) \lellip rx:{\xrad} ry:{\yrad}
        \realadd {15} { \yrad} \outB
        \realadd {15} {-\yrad} \outA
        \realadd {30} { \yrad} \inD
        \realadd {30} {-\yrad} \inC
        \realadd { 0} { \yrad} \inB
        \realadd { 0} {-\yrad} \inA
        \move (0 {\inD}) \clvec (-25 {\inD})(-20 {\outB})(-40 {\outB}) 
        \move (0 {\inB}) \clvec (-18 {\inB})(-18 {\inC})(0 {\inC}) 
        \move (0 {\inA}) \clvec (-25 {\inA})(-20 {\outA})(-40 {\outA}) 
        \savepos (0 0)(*ex *ey)
	  \esegment
	  \move (*ex *ey)

    \bsegment
        \move (40 15) \lellip rx:{\xrad} ry:{\yrad}
        \realadd {15} { \yrad} \outB
        \realadd {15} {-\yrad} \outA
        \realadd {30} { \yrad} \inD
        \realadd {30} {-\yrad} \inC
        \realadd { 0} { \yrad} \inB
        \realadd { 0} {-\yrad} \inA
        \move (0 {\inD}) \clvec (25 {\inD})(20 {\outB})(40 {\outB}) 
        \move (0 {\inB}) \clvec (18 {\inB})(18 {\inC})(0 {\inC}) 
        \move (0 {\inA}) \clvec (25 {\inA})(20 {\outA})(40 {\outA}) 
        \savepos (40 15)(*ex *ey)
    \esegment
    \move (*ex *ey)

	\esegment
    \move (*ex *ey)
}

\newcommand{\idFA}{%
	\bsegment
		\realmult {\fatness} {2.5} \xrad
		\realmult {\fatness} {6.25} \yrad
		\move (0 0) \lellip rx:{\xrad} ry:{\yrad}
		\move (40 0) \lellip rx:{\xrad} ry:{\yrad}
		\realadd { 0} { \yrad} \inB
		\realadd { 0} {-\yrad} \inA
		\move (0 {\inB}) \rlvec (40 0) 
		\move (0 {\inA}) \rlvec (40 0) 
		\savepos (40 0)(*ex *ey)
	\esegment
	\move (*ex *ey)
}

\makeatletter
\renewcommand\section{\@startsection {section}{1}{\z@}%
                                  {-10.5pt \@plus -1pt \@minus -.2pt}%
                                  {3.5pt \@plus.2ex}
                                   {\bf}}
\makeatother

\def\BBdot{ \fcir f:0 r:2 }
\def\BBBdot{ \fcir f:0 r:3 }
\def\BBWdot{ \bsegment \fcir f:1 r:3 \lcir r:3 \esegment }

\newcommand{\R}{\mathbf{R}}
\providecommand{\kat}[1]{\text{\textbf{\textsl{#1}}}}

\newcommand{\Cospan}{\kat{Cospan}}
\newcommand{\Cob}{\kat{2Cob}}

\newcommand{\from}{\leftarrow}

\newcommand{\Addresses}{{
  \bigskip
  \footnotesize

  J. Kock, \textsc{Departament de Matem\`{a}tiques, 
Universitat Aut\`{o}noma de Barcelona, 
SPAIN}\par\nopagebreak
  \textit{E-mail address}: \texttt{kock@mat.uab.cat}

  \medskip

  D.I.~Spivak, \textsc{Department of Mathematics, Massachusetts Institute of Technology,
USA}\par\nopagebreak
  \textit{E-mail address}: \texttt{dspivak@gmail.com}

}}

\begin{document}
\title{
Homotopy composition of cospans
}

\author{
Joachim Kock%
\
and David I. Spivak%
}

\date{}
\maketitle

\vspace*{-1em}

\begin{abstract}
  It is well known that the category of finite sets and cospans, composed by
  pushout, contains the universal {\em special} commutative Frobenius
  algebra.  In this note we observe that the same construction yields also
  general commutative Frobenius algebras, if just the pushouts are changed to
  homotopy pushouts.
\end{abstract}


\section{Introduction}

A Frobenius algebra, or more generally, a Frobenius object in a monoidal
category~\cite{Kock:FA-2DTQFT}, is {\em special} when the comultiplication
followed by the multiplication is the identity map.  
A symmetric monoidal category containing the universal special commutative
Frobenius algebra was exhibited by Lack~\cite{Lack:TAC13}, and independently by
Rosebrugh, Sabadini and Walters~\cite{Rosebrugh-Sabadini-Walters:TAC6} in their
study of general processes: it is the category whose objects are finite sets and whose
morphisms are cospans; the composition of cospans is given by pushout.  Special
commutative Frobenius objects and cospan categories have also been studied
recently by Baez and Fong~\cite{Baez-Fong:1504.05625} in the context of
electrical network theory, and Coya and Fong~\cite{Coya-Fong:1601.02307}
have shown that the 
category of jointly surjective cospans contains the universal {\em extraspecial}
commutative Frobenius algebra, meaning that also unit-followed-by-counit is required to
be the identity.\footnote{In the context of conformal field theory (see for 
example~\cite{Fuchs-Runkel-Schweigert:0110133}),
``extraspecial'' is called just ``special''.}

In this note we observe that ``specialness'' is a result of the simple-minded nature
of pushouts in the category of sets: if one replaces
pushouts by homotopy pushouts, the resulting cospan category
contains the universal commutative Frobenius object, rather than the special one.
To make sense of the homotopy pushout, a homotopical
setting is required.  The minimal setting in which the construction works is that of finite
$1$-dimensional CW complexes up to homotopy, i.e.~topological graphs, the
endpoints of the cospans being still just finite sets.  Graph cospans like this
have a long history in the theory of graph rewriting (an early explicit
reference being~\cite{Gadducci-Heckel}), but hitherto, as far as we know, only
set-theoretical pushouts have been considered.

In light of the growing importance of homotopy viewpoints in
modern mathematical sciences, we wonder whether the finer composition of
cospans, provided by the homotopy pushout, may be of relevance in the various
application areas of the cospan model.
We leave the investigation of these issues to the experts in the respective
fields.

\section{Cobordisms versus cospans}

To put the cospan construction in perspective, and to explain its homotopy
version, we exploit the standard result that the universal commutative Frobenius
algebra lives in the symmetric monoidal category $\kat{2Cob}$ of 2-dimensional cobordisms
\cite[Theorem~3.6.19]{Kock:FA-2DTQFT}.  (A corollary is the notorious equivalence between 2d
topological quantum field theories and commutative Frobenius algebras.)

The
result essentially amounts to writing the axioms for commutative Frobenius 
algebras in graphical language, and observing that they precisely express
evident topological properties of cobordisms.  We shall not repeat
the arguments (nor even the axioms) here, but content ourselves to list the
graphical building blocks corresponding to multiplication and unit, 
comultiplication and counit, and the corresponding cospans:
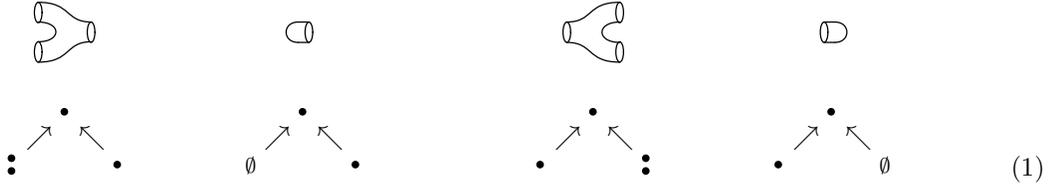
\begin{equation}\label{eqn:four_pics}\begin{texdraw}\setunitscale 0.5 
  \move (20 -15) \multFA 
  \move (225 0) \unitFA 
  
  \move (420 0) \comultFA 
  \move (615 0) \counitFA 

  \move (0 -100)
  \bsegment
  \bsegment
  \move (0 0) \bsegment \move (0 5) \BBBdot \move (0 -5) \BBBdot \esegment
  \move (40 40) \BBBdot
  \move (80 0) \BBBdot
  \htext (20 20) {$\nearrow$}
  \htext (60 20) {$\nwarrow$}
  \esegment
  
  \move (180 0)
  \bsegment
  \move (0 0) \htext{\footnotesize $\emptyset$}
  \move (40 40) \BBBdot
  \move (80 0) \BBBdot
  \htext (20 20) {$\nearrow$}
  \htext (60 20) {$\nwarrow$}
  \esegment
  
  \move (400 0)
  \bsegment
  \move (80 0) \bsegment \move (0 5) \BBBdot \move (0 -5) \BBBdot \esegment
  \move (40 40) \BBBdot
  \move (0 0) \BBBdot
  \htext (20 20) {$\nearrow$}
  \htext (60 20) {$\nwarrow$}
  \esegment
  
  \move (580 0)
  \bsegment
  \move (80 0) \htext{\footnotesize $\emptyset$}
  \move (40 40) \BBBdot
  \move (0 0) \BBBdot
  \htext (20 20) {$\nearrow$}
  \htext (60 20) {$\nwarrow$}
  \esegment
  \esegment
\end{texdraw}\end{equation}
The analogy shown in \eqref{eqn:four_pics} between 2d cobordisms and cospans $A_0\to X\from A_1$ of finite
sets is clear: $A_0$ corresponds to the set of in-boundaries, $A_1$
corresponds to the set of out-boundaries, and $X$ corresponds to the set
of connected components of a cobordism.  Composition of cobordisms is given 
by gluing at the
boundaries (in fact a pushout operation); cospans compose by pushout.

However, the analogy breaks down
quickly in more complicated situations,
because the set $X$ does not contain information 
about the genus of each component of a cobordism. 
The crucial difference can be
pinpointed to the following basic composition of cospans of sets,
expressing precisely the axiom characterizing special Frobenius algebras
among all Frobenius algebras:
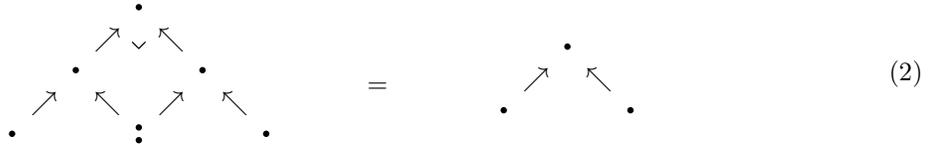
\begin{equation}\label{pushout-set}\parbox{3.3in}{
  \begin{texdraw} \setunitscale 0.6
    \move (0 0) \BBdot
    \move (80 0) \bsegment \move (0 4) \BBdot \move (0 -4) \BBdot \esegment
    \move (160 0) \BBdot
    \move (40 40) \BBdot
    \move (120 40) \BBdot
    
    \htext (20 20){$\nearrow$}
    \htext (60 20){$\nwarrow$}
    \htext (100 20){$\nearrow$}
    \htext (140 20){$\nwarrow$}
    
    \move (80 80) \BBdot
    \htext (60 60){$\nearrow$}
    \htext (100 60){$\nwarrow$}
    
    \move (80 70) \bsegment \move (-4 -12) \lvec (0 -16) \lvec (4 -12) \esegment

    \htext (230 30){$=$}

    \move (310 15) 
    \bsegment
    \move (0 0) \BBdot 
    \move (80 0) \BBdot
    \move (40 40) \BBdot
    \htext (20 20){$\nearrow$}
    \htext (60 20){$\nwarrow$}
    \esegment
  \end{texdraw}}
\end{equation}
We see that
cospans (in the category of sets) cannot render the idea of a ``hole'' (as 
resulting 
from the corresponding composition of cobordisms, 
\raisebox{-5pt}{
\begin{texdraw} \setunitscale 0.35 \move (0 0) \comultFA \multFA \htext (120
0){$=$} \move (160 0) \handleFA \htext (280 0){$\neq$} \move (320 0) \idFA
  \end{texdraw}
}).


\section{Homotopy pushouts and finite $1$-dimensional CW-complexes}

The general idea of homotopy pushouts and homotopy quotients, which has become
an essential slogan in higher category theory~\cite{Baez-Dolan:finset-feynman},
and more recently in type theory~\cite{HoTT-book},
is that instead of equating elements, one sews in a path between them, so as to
keep track of the fact that two elements could be equal in more than one way.
When applied to the composition of cospans, 
the homotopy-pushout operation matches well
with the various interpretations of cospans in applications: to connect two
``devices'', draw
``cables'' from the outputs of the first to the inputs of the second.

Let $X$ and $Y$ be topological spaces, and let $A$ be a finite set, considered 
as a discrete space.
The {\em homotopy pushout} $X \amalg_A Y$ of a span $X \stackrel f \leftarrow A 
\stackrel g \to Y$,
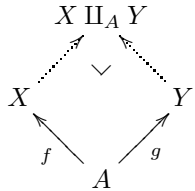
$$\xymatrix @!=7pt {
& X \amalg_A Y \dpullback & \\
   X \ar@{..>}[ru] && Y \ar@{..>}[lu] \\
   &  \ar[lu]^f A \ar[ru]_g &
}$$
is obtained by attaching an interval in $X\amalg Y$ from $f(a)\in X$
to $g(a)\in Y$, for each point $a\in A$.  
For example, the homotopical version of the pushout in \eqref{pushout-set} is shown below, clearly corresponding to \raisebox{-5pt}{\begin{texdraw}\setunitscale 0.35 
\handleFA \end{texdraw}}\,: 
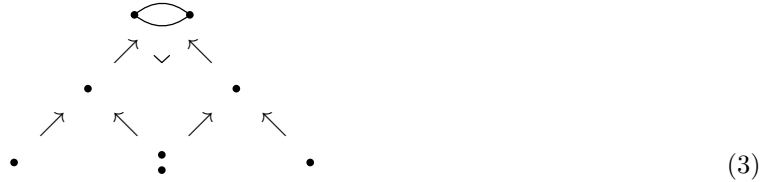
\begin{equation}\label{pushout-top}
  \begin{texdraw} \setunitscale 0.7
    \move (0 0) \BBdot
    \move (80 0) \bsegment \move (0 4) \BBdot \move (0 -4) \BBdot \esegment
    \move (160 0) \BBdot
    \move (40 40) \BBdot
    \move (120 40) \BBdot
    
    \htext (20 20){$\nearrow$}
    \htext (60 20){$\nwarrow$}
    \htext (100 20){$\nearrow$}
    \htext (140 20){$\nwarrow$}
    
    \move (65 80) \BBdot
    \move (95 80) \BBdot
    \move( 65 80) \clvec (75 88)(85 88)(95 80)
    \move( 65 80) \clvec (75 72)(85 72)(95 80)
    
    \htext (60 60){$\nearrow$}
    \htext (100 60){$\nwarrow$}
    
    \move (80 70) \bsegment \move (-4 -12) \lvec (0 -16) \lvec (4 -12) \esegment
  \end{texdraw}
\end{equation}

\vspace*{-4pt}

We see in particular, that even when the spaces $X$, $Y$, and $A$ are all discrete,
the homotopy pushout $X \amalg_A Y$ will not in general be discrete.
To accommodate the homotopy pushout, we need a category
of cospans with discrete endpoints, but with apex 
at least a homotopy $1$-type, say a $1$-dimensional CW complex.
Recall that a {\em finite $1$-dimensional CW-complex}
is a space obtained from a finite set of points by attaching a finite set of
segments by their endpoints.  The homotopy pushout described above is a 
paradigmatic example.
A finite $1$-dimensional CW-complex is thus a topological graph; if it is 
connected it is characterized up to homotopy by
its first Betti number, the number of independent cycles. Any connected
finite $1$-dimensional CW-complex is therefore homotopy equivalent to a bouquet
of circles. Note that the homotopy pushout is a homotopy-invariant construction: replacing
any of the three objects in the span $X\from A\to Y$ by something homotopy
equivalent will not change the homotopy class of the homotopy pushout.

Let $\Cospan$ denote the symmetric monoidal category whose objects are finite sets, and whose
arrows are homotopy classes of cospans $A_0 \to X \leftarrow A_1$, where $X$ is a
finite $1$-dimensional CW-complex.  
The monoidal structure is given by disjoint union.
It is easy to see that such cospans are closed under composition given by
homotopy pushout, since this operation attaches $1$-cells, but nothing
higher-dimensional.  Since we only consider cospans up to homotopy, it is also
clear that composition is associative and unital.

\section{Homotopy cospans and cobordisms}

\medskip
\noindent {\bf Theorem.} {\em
The symmetric monoidal category $\Cospan$, whose objects are 
finite sets and whose morphisms are homotopy classes of
cospans, is the free symmetric monoidal category
on a commutative Frobenius object.
}
\medskip

The universal property can be
established directly, without reference to cobordisms,
by copying over the proof in 
\cite{Kock:FA-2DTQFT} almost verbatim.  This is possible since
the generators and relations of $\Cospan$ as a symmetric monoidal category
are readily seen to correspond exactly to the axioms for commutative Frobenius 
algebras.  The key point here is that the invariants detecting whether
two connected cospans are homotopic are the same as those
classifying connected $2$-cobordisms, as used in \cite{Kock:FA-2DTQFT}.

\medskip

It is also possible to establish this result without mention of 
generators and relations, by means of a geometric formalization of the
analogy between cospans and cobordisms, establishing an equivalence of 
symmetric monoidal categories
$\Cospan\simeq\Cob$.  This provides an alternative proof of the Theorem, since
$\Cob$ is already known to have the universal property (\cite{Kock:FA-2DTQFT},
Theorem~3.6.19).

We first need to enhance the geometric 
realization of a cospan $A_0 \to X \from A_1$, 
by separating out its input and output points.  We do that
systematically by replacing $X$ by the
{\em mapping cylinder} of $A_0+A_1 \to X$; this amounts to attaching a $1$-cell
from each point in $A_0+A_1$ to its image point in $X$, and produces
thus a homotopy equivalent graph whose input and output points are all pairwise disjoint.

As an example, here is the replacement for the
cospan representing multiplication:
\begin{equation}\label{eqn:mappingcylinder}\begin{texdraw}\setunitscale 0.6 
  \move (0 0)
  \bsegment
  \move (0 0) \bsegment \move (0 5) \BBWdot \move (0 -5) \BBWdot \esegment
  \move (40 40) \BBBdot
  \move (80 0) \BBWdot
  \htext (20 20) {$\nearrow$}
  \htext (60 20) {$\nwarrow$}
  \esegment

  \htext(155 20){$\leadsto$}
  \move (240 0)
  \bsegment
  \move (0 0) \bsegment \move (0 5) \BBWdot \move (0 -5) \BBWdot \esegment
  \move (55 50) 
  \bsegment 
  \move (0 0) \BBBdot
  \lvec (-14 7) \move (-16 8) \BBWdot 
  \move (0 0) \lvec (-14 -7) \move (-16 -8) \BBWdot
  \move (0 0) \lvec (13 0) \move (16 0) \BBWdot 
  \esegment
  \move (110 0) \BBWdot
  \htext (20 20) {$\nearrow$}
  \htext (90 20) {$\nwarrow$}
  \esegment
\end{texdraw}\end{equation}
already considerably strengthening the analogy with
\raisebox{-5pt}{
\begin{texdraw} \setunitscale 0.35 \move (0 0) \multFA \end{texdraw}
  } .
(At this point it is worth recalling that a cobordism is also in fact a kind of cospan 
(cf.~\cite{Kock:FA-2DTQFT}, \S 1.2): 
it is the surface with boundary, {\em together with} specific
inclusions $\Sigma_0 \to M \leftarrow \Sigma_1$ from closed $1$-manifolds onto
the boundary.)

With this preparatory mapping-cylinder step, we 
can describe the correspondence geometrically in Euclidean $3$-space.
Given a cospan we choose a suitable embedding of the mapping-cylinder graph
into $\R^2 \times [0,1] \subset \R^3$, with
its input points at level $0$ and output points at level $1$.  
By suitable, we mean: smooth on each $1$-cell, and at each vertex
 the adjacent $1$-cells should not share tangent directions. Other than these requirements, it is not 
important how exactly the embedding is arranged---even the isotopy
class does not matter. In order to construct a cobordism from this, we first take the boundary of a
tubular neighborhood of the embedded graph as a subset of $\R^3$; specifically we can take
the set of points in $\R^3$ at fixed distance $\epsilon>0$ from the graph. If 
$\epsilon$ is chosen small enough, which is possible by compactness of the finitely many closed cells,
the resulting subset will be a topological surface. We intersect this surface with $\R^2\times[0,1]$
to obtain a cobordism, whose boundary consists of small circles around the input
and output points of the graph.

Conversely, any cobordism can be embedded in $\R^2 \times [0,1]$
with in-boundaries at level $0$ and out-boundaries at level $1$. Its
interior admits a deformation retract onto a graph, 
commuting (at the boundary levels $\R^2\times\{0,1\}$) with the deformation
retract of the open disks onto their center points;
the result is unique up to homotopy. With boundary points determined as those at
levels $0$ and $1$, this defines a cospan.
It is clear that these two constructions are mutually
inverse (up to the equivalences in question: homotopy equivalence
of graph cospans, and homeomorphism rel the boundary for cobordisms).

Finally, it is clear from the geometry that composition of cospans corresponds
precisely to composition of cobordisms: Given composable cospans $A_0\to X\from
A_1$ and $A_1\to Y\from A_2$, the composite given by the homotopy
pushout can be realized geometrically by embedding (the mapping cylinders of) $X$
into $\R^2\times [0,1]$ and $Y$ into $\R^2\times [2,3]$ and sewing in $1$-cells
$A_1\times I$ in $\R^2\times [1,2]$.  Cobordisms can be composed in exactly the
same way, by sewing in cylinders in $\R^2\times[1,2]$ between the
boundary circles at level $1$~and~$2$.

\bigskip
\noindent
{\bf Acknowledgments.} Kock was supported by grant MTM2013-42293-P of Spain, and Spivak was supported by AFOSR grant FA9550--14--1--0031 and NASA grant NNH13ZEA001N. We thank the referee for useful remarks.

\Addresses
\vspace{-.2in}

\end{document}